\documentclass{amsart}

\usepackage{amsthm,amssymb,latexsym,amsmath}
\usepackage{mathrsfs}
\usepackage{graphicx}
\input xy
\xyoption{all}
\usepackage[all]{xy}

\renewcommand{\dim}{\mathrm{dim}}
\newcommand{\sF}{ {\mathscr F}}

\newcommand{\Q}{\mbox{${\mathcal Q}$}}

\newcommand{\G}{\mbox{${\mathscr G}$}}

\def\H{\mathcal H}
\renewcommand{\H}{\text{\rm H}}

\newtheorem{nada}{Nada}[section]

\newtheorem{thm}[nada]{Theorem}
\newtheorem{lemma}[nada]{Lemma}

\newtheorem{example}[nada]{Example}

\newcommand{\bc}{\begin{center}}
\newcommand{\ec}{\end{center}}
\newcommand{\noi}{\noindent}

\begin{document}

\title{Determination of Baum-Bott residues of  Higher codimensional foliations }
\hyphenation{ho-mo-lo-gi-cal}
\hyphenation{fo-lia-tion}

\begin{abstract}
Let $\sF$ be a singular holomorphic foliation, of codimension $k$, on a complex compact  manifold such that its singular set has codimension $\geq k+1$.
In this work we determinate   Baum-Bott residues for  $\sF$  with respect to  homogeneous symmetric polynomials of degree $k+1$.
We drop the  Baum-Bott's generic hypothesis and we show that the residues can be expressed in terms of  the  Grothendieck residue  of an one-dimensional foliation on a $(k+1)$-dimensional  disc transversal to a $(k+1)$-codimensional component of the singular set of $\sF$. Also, we show that  Cenkl's algorithm for non-expected dimensional singularities holds dropping the   Cenkl's regularity assumption.  
\end{abstract}

\author{Maur\'icio Corr\^ea}
\address{Maur\'icio Corr\^ea \\ ICEx - UFMG \\
Departamento de Matem\'atica \\
Av. Ant\^onio Carlos 6627 \\
30123-970 Belo Horizonte MG, Brazil } \email{mauriciojr@ufmg.br}

\author{Fernando Louren\c co}
\address{ Fernando Lourenço \\ DEX - UFLA \\ Av. Doutor Sylvio Menicucci, 1001, Kennedy, 37200000, Lavras, Brasil
}
\email{fernando.lourenco@dex.ufla.br}
\thanks{ }

\maketitle

\section{Introduction}

In \cite{BB}  P. Baum and R. Bott  developed  a   general residue theory
for singular holomorphic foliations on complex manifolds. More precisely, they proved the following result:

\begin{thm}[Baum-Bott] Let $\sF$ be a holomorphic foliation of codimension $k$ on a complex manifold $M$ and $\varphi$ be a homogeneous symmetric polynomials
of degree $d$ satisfying $ k < d \leq n$. Let $Z$ be a compact connected component of the singular set $\mathrm{Sing}(\sF)$. Then,  there exists a homology class $\mathrm{Res}_{\varphi}(\sF, Z) \in \H_{2(n - d)}(Z; \mathbb{C})$ such that:
\begin{enumerate}
\item[$i)$] $\mathrm{Res}_{\varphi}(\sF, Z) $ depends only on $\varphi$ and on the local behavior of the leaves of $\sF$ near $Z$,
\item[$ii)$] Suppose that  $M$ is compact and denote by $\mathrm{Res}({\varphi}, \sF, Z):=\alpha_{\ast}\mathrm{Res}_{\varphi}(\sF, Z)$,
where $\alpha_{\ast}$ is the composition of the maps
$$\displaystyle \H_{2(n - d )}(Z; \mathbb{C}) \stackrel{i^{\ast}}\longrightarrow \H_{2(n - d )}(M; \mathbb{C}) $$ and $$\displaystyle \H_{2(n - d )}(M; \mathbb{C}) \stackrel{P}\longrightarrow \H^{2d}(M; \mathbb{C}) $$
 with $i^{\ast}$ is the induced map of inclusion $i : Z \longrightarrow M$ and $P$ is the Poincar\'e duality. Then 
 $$\displaystyle \varphi (\mathcal{N}_{\sF}) = \sum_{Z} \mathrm{Res}({\varphi}, \sF, Z).$$
 \end{enumerate}

\end{thm}

The computation and determination  of the residues is   difficult in general.  
If the foliation $\sF$ has  dimension one with isolated singularities ,  Baum and Bott  in \cite{BB1}  show that residues can be expressed in terms of a  Grothendieck residue, i.e, for each $p\in \mathrm{Sing}(\sF)$ we have 
$$
\mathrm{Res}_{\varphi}(\sF, Z)=  \mbox{Res}_{p}\Big[\varphi(JX)\frac{dz_{1}\wedge\cdots   \wedge dz_{n}}{X_{1}\cdots   X_{n}}\Big] ,
$$
where $X$ is a germ  of holomorphic vector field at $p$ tangent to $\sF$ and $JX$ is the jacobian of $X$.

The subset of $\mathrm{Sing}(\sF)$ composed by analytic subsets of codimension $k+1$  will be denoted by   $\mathrm{Sing}_{k+1}(\sF)$ and it is called     {\it  the  singular set of $\sF$ with expected codimension }.  Baum and Bott  in \cite{BB} exibes the residues for generic componentes of Sing$_{k + 1}(\sF)$. Let us recall this result:

An  irreducible component $Z$ of $\mathrm{Sing}_{k + 1}(\sF)$ comes endowed with a filtration. For given point  $p \in Z$ 
choose holomorphic vector fields $v_{1},\dots, v_{s}$ defined on an open neighborhood $U_{p}$ of $p \in M$ and such that for all $x \in U_{p}$, the germs at $x$ of the holomorphic vector fields $v_{1},\dots,v_{s}$ are in $\sF_{x}$ and span $\sF_{x}$ as a $\mathcal{O}_{x}$-module. Define a subspace $V_{p}(\sF) \subset T_{p}M$ by letting $V_{p}(\sF)$ be the subspace of $T_{p}M$ spanned by
 $v_1(p),\dots v_s(p)$. We have
$$ Z^{(i)} = \{ p \in Z ; \ \ \dim (V_{p}(\sF)) \leq n - k - i \} \ \ \mbox{for} \ \ i = 1,\dots,n - k. $$
 Then,
$$ Z \supseteq Z^{(1)} \supseteq Z^{(2)} \supseteq \dots \supseteq Z^{(n - k)}$$
is a filtration of $Z$.
Now, consider   a symmetric homogeneous polynomial  $\varphi$   of degree $k+1$. 
Let  $Z \subset \mathrm{Sing}_{k+1}(\sF)$ be an irreducible component.  Take a generic point $p \in 
Z$ such that $p$ is a point where $Z$ is smooth and disjoint from the other singular
components. Now, consider  $B_{p}$ a ball centered at $p $,  of dimension $k+1$ sufficiently small and transversal to $Z$ in $p$.
In \cite[Theorem 3, pg 285]{BB}    Baum and Bott proved   under the following generic  assumption 
$$ \mbox{cod}(Z) = k + 1 \ \ \ \mbox{and} \ \ \  \mbox{cod}(Z^{(2)} ) < k + 1 $$
that we have 
$$ \mathrm{Res}(\sF, \varphi ; Z) = \mathrm{Res}_{\varphi}(\sF|_{B_p}; p)[Z],  $$
where $ \mathrm{Res}_{\varphi}(\sF|_{B_p}; p)$ represents  the Grothendieck residue at $p$ of the one dimensional  foliation $\sF|_{B_p}$ on $B_p$ and $[Z]$ denotes the integration current associated to $Z$. 

In \cite{BruPer} and \cite{CoFer} the authors  determine the residue $\mathrm{Res}(\sF, c_1^{k+1} ; Z)$, but even in this case they do not show that we can calculate these residues in terms of  the Grothendieck residue  of a foliation on a  transversal disc. In \cite{Vis} Vishik proved the same result 
under the  Baum-Bott's generic hypotheses but supposing that the foliation has locally free tangent sheaf.  In \cite{BS} F. Bracci and T. Suwa    study the behavior of the Baum-Bott residues under smooth
deformations, providing  an effective way of computing residues.

In this work we drop the Baum-Bott's generic hypotheses  and we prove the following  :

\begin{thm} \label{principal} Let $\sF$ be a singular holomorphic foliation of codimension $k$ on a compact complex manifold $M$ such that  $\mathrm{cod} (\mathrm{Sing}(\sF) )\geq k + 1$. Then,
$$ \mathrm{Res}(\sF, \varphi ; Z) = \mathrm{Res}_{\varphi}(\sF|_{B_p}; p)[Z], $$
 where  $\mathrm{Res}_{\varphi}(\sF|_{B_p}; p)$ represents  the Grothendieck residue at $p$ of the one dimensional  foliation $\sF|_{B_p}$ on a $(k+1)$-dimensional  transversal ball $B_p$.
\end{thm}

Finally, in the last  section we apply   Cenkl's algorithm for non-expected dimensional singularities \cite{Ce}. Moreover, we  drop   Cenkl's  regularity  hypothesis and we conclude  that  it is  possible to calculate the residues for foliations whenever
  $\mathrm{cod} (\mathrm{Sing}(\sF)) \geq k+s$, with $s\geq 1$. 
  
  \subsection*{Acknowledgments}
We are grateful to   Jean-Paul Brasselet, Tatsuo Suwa and  Marcio G.  Soares    for interesting conversations.
This work was partially supported by CNPq, CAPES, FAPEMIG and FAPESP-2015/20841-5.
We are grateful to   Institut de Math\'ematiques de Luminy- Marseille   and Imecc--Unicamp for   hospitality. 
Finally, we would like to thank the referee by the suggestions, comments
and improvements to the exposition. 

\section{Holomorphic foliations}

Denote by $ \Theta_M$ the tangent sheaf of $M$. 
A foliation $\sF$ of  codimension $k$ on an $n$-dimensional complex manifold $M$ is given by a exact sequence of coherent  sheaves 
$$
  0  \longrightarrow T\sF  \longrightarrow  \Theta_M \to N_{\sF} \longrightarrow 0,
$$
such that  $ [T\sF,T\sF] \subset  T\sF$ and the normal sheaf $N_{\sF}$ of $\sF$ is a   torsion free sheaf of   rank $k\leq n -1$.
The sheaf $T\sF$ is called the tangent sheaf of $\sF$. The singular set of $\mathscr{F}$ is defined by $\mathrm{Sing}(\sF):=\mathrm{Sing}(N_{\sF})$. 
The dimension of $\sF$ is $\dim(\sF)=n-k$.

Also, a foliation $\sF$,  of  codimension $k$,  can be induced by a exact sequence
$$
  0  \longrightarrow N_{\sF}^{\vee}    \longrightarrow  \Omega_M^1\to \Q_{\sF}  \longrightarrow 0,
$$
 where $\Q_{\sF}$ is a torsion free sheaf of   rank $n-k$. Moreover,  the singular set of $\mathscr{F}$ is  $\mathrm{Sing}(\Q_{\sF})$. Now, by taking the wedge product of the map $  N_{\sF}^{\vee}    \longrightarrow  \Omega_M^1 $ we get a morphism
 $$
\bigwedge^k  N_{\sF}^{\vee}    \longrightarrow  \Omega_M^k
 $$
and twisting by $(\bigwedge^k  N_{\sF}^{\vee})^{\vee}=\det(N_{\sF})$ we obtain a morphism
 $$
\omega: \mathcal{O}_M \longrightarrow  \Omega_M^k\otimes  \det(N_{\sF}).
 $$
Therefore, a foliation is induced by a twisted holomorphic  $k$-form
$$\omega \in \H^0(X,\Omega_M^k\otimes \det(N_{\sF}) )$$ which is locally decomposable outside the singular set of $\sF$.
That is, by the classical Frobenius Theorem for each point $p\in X\setminus  \mathrm{Sing}(\sF)$ there exists a neighbourhood $U$ and
holomorphics $1$-forms $\omega_1, \dots,\omega_k \in \H^0(U, \Omega_U^1)$ such that
$$
\omega|_{U}=\omega_1 \wedge \cdots   \wedge \omega_k
$$
and
$$
d \omega_i \wedge \omega_1 \wedge \cdots   \wedge \omega_k=0
$$
for all $i=1,\dots,k$.

\section{Proof of  the Theorem}

Given a multi-index $\alpha = (\alpha_{1},\dots,\alpha_{k})$ with $\alpha_{j} \geq 0 $ for $j = 1,\dots,k$,  consider  the homogeneous symmetric polynomial of degree $k + 1$, $\varphi = c_{1}^{\alpha_{1}}c_{2}^{\alpha_{2}}\cdots   c_{k}^{\alpha_{k}} $  such that $1\alpha_{1} + 2\alpha_{2} + \cdots    + k\alpha_{k} = k + 1 $. 

 Let us consider the  twisted $k$-form $ \omega \in  \H^{0}(M, \Omega^{k}_{M} \otimes \det(N_{\sF}))$ induced by $\sF$.
Denote by Sing$_{ k+1}(\sF)$ the union of the irreducible components of $\mathrm{Sing}(\sF)$ of pure codimension $k + 1$. 
Consider an open subset $U\subset M\setminus \mathrm{Sing}(\sF) $. Thus, the form $\omega|_U$ is decomposable and integrable. That is,   $\omega|_U$ is given by a product of $k$ 1-forms $\omega_{1} \wedge \cdots    \wedge \omega_{k}$. Then, it is possible to find a matrix of $(1,0)$-forms $(\theta_{ls}^*)$ such that
$$ \partial \omega_{l} = \sum_{s = 1}^{k} \theta_{ls}^* \wedge \omega_{s}, \ \   \overline{\partial} \omega_{l} =0, \ \ \forall \ \  l = 1,\dots,k.$$
We have that $\omega_{1} , \ldots   ,  \omega_{k}$ is a local frame for $N_{\sF}^*|_{U}$ and the identity above induces on $U$ the {\it  Bott partial connection}
$$
\nabla : C^{\infty }( N_{\sF}^*|_{U} ) \to  C^{\infty }( (T\sF^*\oplus \overline{TM}) \otimes N_{\sF}^* |_{U} )
$$
defined by
$$
\nabla_v(\omega_{l})= i_{v} ( \partial \omega_{l}), \ \   \ \  \nabla_u(\omega_{l})= i_{u} ( \overline{\partial} \omega_{l})=0,
$$
where $v \in C^{\infty }( T\sF|_{U} )$ and $u\in C^{\infty }(   \overline{TM}  |_{U} )$ which can be    extended  to a connection 
$D^*:C^{\infty }( N_{\sF}^*|_{U} ) \to  C^{\infty }( (TM^*\oplus \overline{TM}) \otimes N_{\sF}^* |_{U} )$ in the following way
$$
D^*_v(\omega_{l})=\sum_{s=1}^{k} i_v(\pi(\theta_{ls}^*))  \omega_{s} , \ \ \ \ D^*_u(\omega_{l})= i_{u} ( \overline{\partial} \omega_{l})=0
$$
where $v \in C^{\infty }( TM|_{U} )$ and $u\in C^{\infty }(   \overline{TM}  |_{U} )$ and $\pi:TM^*|_{U}\to N_{\sF}^*|_{U} $ is the natural projection. 
Let $\theta^{*} $ be the matrix of the  connection $D^*$, then $\theta:=[-\theta^{*}]^t$ is 
the matrix of the  induced  connection $D$  with respect to the
frame $\{ \omega_{1} ,  \dots   , \omega_{k}\}$.

Let $K$ be the curvature of the   connection $D$ of  $N_{\sF}$ on $M\setminus \mathrm{Sing}(\sF)$.  
It follows from Bott's vanishing Theorem \cite[Theorem 9.11, pg 76]{Suwa} that  $\varphi(K) =0 $ .  Let $V$ be  a small neighborhood of  $\mathrm{Sing}_{k + 1}(\sF)$.
 We regularize  $\theta$ and $K$ on $V$ , i.e. we choose a matrix of 
smooth  forms $\widehat{\theta}$ and $\widehat{K}$ coinciding with $\theta$ and $K$ outside of  $V$, respectively.
By hypothesis   $\dim(\mathrm{Sing}(\sF))\leq n- k-1$ we conclude by a dimensional reason  that, for $\deg(\varphi) =  k+1$, only the components of dimension $n-k -1$ of $\mathrm{Sing}(\sF)$ 
 play a role. In fact,  since 
$\mathrm{Res}_{\varphi}( \sF, Z) \in  \H_{2(n-k-1)}(Z, \mathbb{C})$, components of dimension smaller than $n-k -1$ contribute nothing. 
This means that  $\varphi(\widehat{K})$ localizes on  Sing$_{k + 1}(\sF)$. Then, $\varphi(\widehat{K})$  has compact support on  $V$, where $V$ is a small neighborhood of 
$\mathrm{Sing}_{k + 1}(\sF)$.
That is,
$$
\mathrm{Supp}(\varphi(\widehat{K})) \subset \overline{V} .
$$
Then
$$
\varphi(\widehat{K})=  \sum_{Z_i}  \widehat{\lambda_i}(\varphi) [Z_i],
$$
where $Z_i$ is an irreducible component of Sing$_{k+1}(\sF)$ and $\widehat{\lambda_i}(\varphi) \in \mathbb{C}$. On the other hand, we have that 
$$
\varphi(N_{\sF})=\sum_{Z_i} \mathrm{Res}({\varphi}, \sF, Z_i) =  \sum_{Z_i} \lambda_i(\varphi) [Z_i].
$$
 We will show that $\lambda_i(\varphi)=\widehat{\lambda_i}(\varphi)$, for all $i$. 
In particular, this implies that   $\varphi(\widehat{K})=\varphi(N_{\sF})$. Thereafter, we will determinate the numbers $\widehat{\lambda_i}(\varphi)$.

Consider the unique complete polarization of the polynomial $\varphi$, denoted by $\widetilde{\varphi}$. That is, $\widetilde{\varphi}$ is a symmetric $k$-linear  function such that 
$$\left(\frac{1}{2 \pi i}\right)^{k+1}  \widetilde{\varphi}(\widehat{K},\dots,\widehat{K}) = \left(\frac{1}{2 \pi i}\right)^{k+1}  \varphi(\widehat{K}). $$

Take a generic point $p \in Z_i$, that is, $p$ is a point where $Z_i$ is smooth and disjoint from the other  components.    Let us consider $L \subset M$ a $(k+1)$-ball intersecting transversally Sing$_{k+1}(\sF)$ at a single
point $p \in Z_i$ and  non intersecting  other  component. 
Define 
\begin{equation}\label{eq.4.1}
BB(\sF,\varphi; Z_i) := \left(\frac{1}{2 \pi i}\right)^{k+1}  \int_{L}     \varphi(  \widehat{K}).
\end{equation}
Then  $\widehat{ \lambda_i}(\varphi) = BB(\sF,\varphi; Z_i)$. In fact
$$
BB(\sF,\varphi; Z_i)=\left(\frac{1}{2 \pi i}\right)^{k+1}  \int_{L} \varphi(\widehat{K})= [L] \cap   [\varphi(  \widehat{K})]= \widehat{ \lambda_i}(\varphi)  [L] \cap [Z_i]=\widehat{ \lambda_i}(\varphi) 
$$
since $[L]  \cap [Z_i]=1$ and $[L]  \cap [Z_i]=0$ for all $i\neq j$.
For each   $j = 1,\dots,k$, define the polynomial  
$$ \varphi_{j}(\widehat{\theta}, \widehat{K}) := \widetilde{\varphi}( \widehat{\theta},\underbrace{-2 \widehat{\theta} \wedge \widehat{\theta} ,\dots,-2\widehat{\theta}\wedge \widehat{\theta}}_{j - 1},\underbrace{ \widehat{K},\dots,\widehat{K}}_{k - j}). $$
Now, we consider the $(2k + 1)$- form
$$\displaystyle \varphi_{\alpha}(\widehat{\theta}, \widehat{K}) = \sum_{j = 0}^{k - 1}(-1)^{j} \frac{(k - 1)!}{2^{j}(k - j - 1)! (k + j)! } \varphi_{j + 1}(\widehat{\theta}, \widehat{K}) .$$
It follows from  \cite[Lemma 2.3, pg 5] {Vis} that  on $X\setminus Sing_{k+1}(\sF)$ we have 
$$ d(\varphi_{\alpha}(\widehat{\theta}, \widehat{K}) ) = \varphi(\widehat{K}). $$
Consider $i: B \to M$ an embedding   transversal to $Z_i$ on $p$ as above, i.e, $i(B)=L$. 
We have then an one-dimensional foliation $\sF|_{L}=i^*\sF$ on $B$ singular only on $i^{-1}(p)=0$. 
We have that 
$$
\widehat{\lambda_i}(\varphi)= BB(\sF,\varphi; Z_i)= \left(\frac{1}{2 \pi i}\right)^{k+1}  \int_{L}    \varphi(  \widehat{K})= \left(\frac{1}{2 \pi i}\right)^{k+1}  \int_{B}    \varphi(i^*  \widehat{K}).
$$
Now, by Stokes's theorem we obtain
$$
\widehat{\lambda_i}(\varphi) =  \left(\frac{1}{2 \pi i}\right)^{k+1}  \int_{B}    \varphi(i^* \widehat{K})= \left(\frac{1}{2 \pi i}\right)^{k+1}  \int_{B}   d(\varphi_{\alpha}(i^*\widehat{\theta},i^* \widehat{K}) )
=\left(\frac{1}{2 \pi i}\right)^{k+1}  \int_{\partial B}    \varphi_{\alpha}(i^*\widehat{\theta},i^* \widehat{K}).
$$
Firstly, it follows from \cite[Lemma 4.6]{Vis}  that
\begin{equation}\label{V}
\widehat{\lambda_i}(\varphi)=\left(\frac{1}{2 \pi i}\right)^{k+1}  \int_{\partial B}    \varphi_{\alpha}(i^*\widehat{\theta},i^* \widehat{K})=  \mbox{ Res}_{\varphi}(i^*\sF ; 0)
\end{equation}
Now, we will adopt the Baum and Bott construction \cite{BB}.  Denote by $\mathcal{A}_M$ the sheaf of germs
of real-analytic functions on $M$. Consider on $M$ a locally free resolution  of  $N_{\sF}$
$$
0\to \mathcal{E}_r  \to  \mathcal{E}_{r-1}  \to \cdots \to  \mathcal{E}_{0} \to N_{\sF}\otimes \mathcal{A}_M\to 0.
$$
Let   $D_{q},D_{q-1}, \dots , D_0$ be  connections for  $\mathcal{E}_q ,  \mathcal{E}_{q-1} , \dots,  \mathcal{E}_{0}$, 
respectively. Set  the curvature of $D_{i}$ by $K_i=K(D_{i})$.  By   using   Baum-Bott notation \cite[pg 297]{BB} we have that 
$$
\varphi(K_{q}| K_{q-1}| \cdots | K_0)= \varphi( N_{\sF}).
$$

Consider on $V$ a locally free resolution of the tangent sheaf of $\sF$:
\begin{equation}\label{seqexa0}
0\to \mathcal{E}_q  \to  \mathcal{E}_{q-1}  \to \cdots \to  \mathcal{E}_{1} \to T\sF \otimes \mathcal{A}_V\to 0.
\end{equation}
Combining this sequence with the sequence 
$$
0\to T\sF \otimes \mathcal{A}_V  \to TV \to N_{\sF}\otimes \mathcal{A}_V\to 0.
$$
we get 
\begin{equation} \label{seqexa2}
0\to \mathcal{E}_q  \to  \mathcal{E}_{q-1}  \to \cdots \to \mathcal{E}_{1} \to TV \to N_{\sF}\otimes \mathcal{A}_V\to 0.
\end{equation}
Pulling back   the sequence (\ref{seqexa0})  by  $i:B\to V$  we  obtain an exact  sequence  on $B$: 
\begin{equation}\label{seqexa}
0\to i^*\mathcal{E}_q  \to i^* \mathcal{E}_{q-1}  \to \cdots \to i^*\mathcal{E}_{1}   \to   i^*(T\sF \otimes \mathcal{A}_V)\to 0.
\end{equation}
Since $B$ is a small ball we have the splitting   $i^*TV=TB\oplus N_{B|V}$,  where   $N_{B|V}$ denotes its  normal bundle. We consider the  projection $\xi:i^*TV \to TB$ and 
we map $i^*TV$ to $N_{i^*\sF}$ via 
$$
i^*TV \stackrel{\xi}{\to} TB \to  N_{i^*\sF}
$$
which give us an exact sequence 
\begin{equation}\label{projec}
 0\to i^*(T\sF \otimes \mathcal{A}_V) \to i^*TV  \to  N_{i^*\sF} \otimes \mathcal{A}_B\to 0
\end{equation}
Now,   combining  the  exact sequences (\ref{seqexa}) and  (\ref{projec})  we obtain an exact sequence 
$$
0\to i^*\mathcal{E}_q  \to i^* \mathcal{E}_{q-1}  \to \cdots \to i^*\mathcal{E}_{1}   \to i^*TV \to N_{i^*\sF}\otimes \mathcal{A}_B\to 0.
$$
Let   $D_{q},D_{q-1}, \dots , D_0$ be  connections for  $\mathcal{E}_q ,  \mathcal{E}_{q-1} , \dots,  \mathcal{E}_{1} , TV$, 
respectively. 
Observe that 
$$
i^*\varphi(K_{q}| K_{q-1}| \cdots | K_0)=\varphi(i^*K_{q}| i^*K_{q-1}| \cdots | i^*K_0)= \varphi( N_{i^*\sF}).
$$
Finally, it follows from \cite[ Lemma 7.16 ]{BB} 
$$
\mbox{ Res}_{\varphi}(i^*\sF ; 0)= \left(\frac{1}{2 \pi i}\right)^{k+1}  \int_{ B}  \varphi(i^*K_{q}| i^*K_{q-1}| \cdots | i^*K_0)=\left(\frac{1}{2 \pi i}\right)^{k+1}  \int_{ B}    i^*\varphi(K_{q}| K_{q-1}| \cdots | K_0)
$$ 
and  \cite[9.12, pg 326]{BB}  that
\begin{equation}\label{lambdaBB}
\mbox{ Res}_{\varphi}(i^*\sF ; 0)= \left(\frac{1}{2 \pi i}\right)^{k+1}  \int_{ B}    i^*\varphi(K_{q}| K_{q-1}| \cdots | K_0) = \lambda_i(\varphi).
\end{equation}


Thus, we conclude from (\ref{V}) and  (\ref{lambdaBB})   that  $\lambda_i(\varphi)=\widehat{\lambda_i}(\varphi),$ for all $i$. This implies that   $\varphi(\widehat{K})=\varphi(N_{\sF})$.

Now, we will determinate the numbers $\widehat{\lambda_i}(\varphi)$.
Let $X=\sum_{r=1}^{k+1} X_i\partial/ \partial z_i$ be a  vector field  inducing $i^*\sF$ on $B$ and $J(X)$ denotes the Jacobian of $X$. Let  $\omega$ be the $1$-form on $B\setminus \{0\}$ such that $i_X(\omega)=1$.  It follows  from \cite[Corollary 4.7]{Vis} that 
$$
\widehat{\lambda_i}(\varphi) = \left(\frac{1}{2 \pi i}\right)^{k+1}  \int_{\partial B}    \varphi_{\alpha}(i^*\widehat{\theta},i^* \widehat{K})= 
 \left(\frac{1}{2 \pi i}\right)^{k+1}  \int_{\partial B}    \omega \wedge (\overline{\partial}\,\omega)^{k}\varphi(-J(X)).
$$
Thus,
$$
\widehat{\lambda_i}(\varphi)=
 \left(\frac{1}{2 \pi i}\right)^{k+1}  \int_{\partial B} (-1)^{k+1}   \omega \wedge (\overline{\partial}\,\omega)^{k}\varphi(J(X)),
$$
By using   Martinelli's formula  \cite[pg. 655]{gri} we have
$$
\widehat{\lambda_i}(\varphi) = 
 \left(\frac{1}{2 \pi i}\right)^{k+1}  \int_{\partial B}   (-1)^{k+1} \omega \wedge (\overline{\partial}\,\omega)^{k}\varphi(J(X))= 
 \mbox{Res}_{0}\Big[\varphi(JX)\frac{dz_{1}\wedge\cdots   \wedge dz_{k+1}}{X_{1}\cdots   X_{k+1}}\Big]. 
$$
Therefore,
$$
\widehat{\lambda_i}(\varphi) =  \mbox{ Res}_{\varphi}(i^*\sF ; 0) = \mbox{Res}_{\varphi}(\sF|_{L}; p),
$$
 where  $\mathrm{Res}_{\varphi}(\sF|_{L}; p)$ represents  the Grothendieck residue at $p$ of the one dimensional  foliation $\sF|_{L}$ on a $(k+1)$-dimensional  transversal ball $L$.

\section{Examples }

In the next examples, with a slight abuse of
notation, we write $ \mathrm{Res}(\sF, \varphi; Z_{i}) = \lambda_i(\varphi). $

\begin{example}\label{ex.4.0.5} Let $\sF$ be the logarithmic foliation on $\mathbb{P}^{3}$ induced, locally in \\ $(\mathbb{C}^{3}, (x,y,z))$ by the polynomial 1-form
$$ \omega = yzdx + xzdy + xy dz. $$
In this chart, the singular set of $\omega$ is the union of the lines   $ Z_{1} = \{x = y = 0   \}; \ \ Z_{2} = \{x = z = 0   \} \ \ \mbox{and} \ \ Z_{3} = \{y = z = 0   \}$. 
We have $  \mathrm{Res}(\sF,c_{1}^{2}; Z_{i}) = \mathrm{Res}_{c_{1}^{2}} (\G; p_{i})$, where $\G$ is a foliation on $D_{i}$ with $D_{i}$ a 2-disc cutting transversally $Z_{i}$.
Consider $D_{1}=\{ ||(x,y) ||\leq 1, \ \  z = 1 \}$ then, we have
$$ \omega|_{D_{1}} =: \omega_{1} = y dx + x dy \ \ \ \ \mbox{with} \ \ \mbox{dual} \ \ \mbox{vector} \ \  \mbox{field} \ \ \ \ X_{1} = x \frac{\partial}{\partial x} - y \frac{\partial}{\partial y}.$$
Then, $D_{1} \cap Z_{1} = \{ p_{1} = (0,0, 1) \}. $ Now, a straightforward calculation shows that

\bc $J X_{1} = \begin{pmatrix} 1  &     0                \\
                           0  & -1  \end{pmatrix}.$ \ec                      
 Thus,
$$ \mathrm{Res}_{c_{1}^{2}}(\G; p_{1}) = \frac{c_{1}^{2}(J X_{1}(p_{1}))}{\det (J X_{1}(p_{1}))} = 0. $$
The same holds for $Z_{2}$ and $Z_{3}$.
The  foliation $\sF$ is induced, in homogeneous coordinates $[X, Y, Z, T]$, by the form
$$ \widetilde{\omega} = YZT dX + XZT dY + XYT dZ -3XYZ dT. $$
The singular set of $\sF$ is the union of the lines $Z_{1}$, $Z_{2}$, $Z_{3}$, and 
\bc
$ Z_{4} = \{ T = X = 0 \}, \ \ Z_{5} = \{ T = Y = 0 \} \ \ \mbox{and} \ \ Z_{6} = \{ T = X = 0 \}. $ 
\ec
For $ Z_{4} = \{ X = T = 0    \}$ we can consider the local chart $U_{y} = \{ Y = 1 \} $. Then, we have,
$$ \omega_{y} := \widetilde{\omega}|_{U_{y}} = z t dx + x t dz - 3xz d t. $$
Take a 2-disc transversal    $D_{2} = \{||(x,t)||\leq 1,\  z = 1 \} $.
$$  \omega_{2} := \omega_{y}|_{D_{2}} = t dx - 3x dt \ \ \ \ \mbox{with} \ \ \mbox{dual} \ \ \mbox{vector} \ \  \mbox{field} \ \ \ \ \displaystyle X_{2} = -3x \frac{\partial}{ \partial x} - t \frac{\partial}{\partial t}.$$
Thus,  $ Z_{4} \cap D_{2} = \{ (0,1,0) =: p_{4} \} $ and 
\bc $J X_{2}(p_{4}) = \begin{pmatrix} -3  &     0                \\
                           0  & -1  \end{pmatrix}.$ \ec

\noi Therefore, $ \mathrm{Res}\displaystyle _{c_{1}^{2}}(\G; p_{4}) = \frac{c_{1}^{2}(JX_{2})(p_{4})}{\det(JX_{2})(p_{4})} = \frac{16}{3}$.
An analogous calculation shows that
$$
\mathrm{Res}\displaystyle _{c_{1}^{2}}(\G; p_{5}) =\mathrm{Res}\displaystyle _{c_{1}^{2}}(\G; p_{6})   = \frac{16}{3}.
$$
 Now, we will verify the formula  
$$
\displaystyle  c_{1}^{2}(N_{\sF})  = \sum_{i = 1}^{6}\mathrm{Res}(\sF, c_{1}^{2} ; Z_{i})[Z_{i}].
$$
On the one hand,  Since $\det(N_{\sF})=\mathcal{O}_{\mathbb{P}^3}(4)$, then 
$$c_{1}^{2}(N_{\sF}) = c_{1}^{2}(\det(N_{\sF}))=16h^2,$$
  where $h$ represents the hyperplane class.
On the other hand, by the above calculations  and since $[Z_{i}]=h^2$, for all $i$, we have 
$$\sum_{i = 1}^{6} \mathrm{Res}(\sF, c_{1}^{2} ; Z_{i})[Z_{i}]= 0[Z_{1}] + 0[Z_{2}] + 0[Z_{3}] + \frac{16}{3}[Z_{4}] + \frac{16}{3}[Z_{5}] + \frac{16}{3}[Z_{6}]=16h^2. $$

\end{example}

The following example is due to D. Cerveau and A. Lins Neto, see  \cite{CerLins}. It originates from the so-called exceptional component of the space of codimension one holomorphic foliations of degree $2$  of $\mathbb{P}^{n}$. We can simplify the computation as done by M. Soares in \cite{Soa}.

\begin{example} Consider $\sF$ be a holomorphic foliation of codimension one on $\mathbb{P}^{3}$, given locally by the 1-form
$$ \omega = z(2y^{2} - 3x  )dx + z(3z - xy)dy - (xy^{2} -2x^{2} + yz )dz. $$

The singular set of this foliation has one connect component, denoted by $Z$, with 3 irreducible components, given by: \\

\noi 1) the twisted cubic $ \Gamma: \ \ \ y \longmapsto (2/3 y^{2}, y, 2/9 y^{3})$, \\

\noi 2) the quadric $Q: \ \ \ y \longmapsto (y^{2}/2, y, 0) $, \\

\noi 3) the line $L: \ \ \ y \longmapsto (0, y, 0). $ \\

We consider a   transversal   $2$-disc   $ D \subset \{ y = 1 \} $ and we take the restriction of $\sF$ on the afine open  $\{ y = 1 \} $. We have an one-dimensional holomorphic foliation, denoted by $\G$, given by the 1-form on $H$
$$ \widetilde{\omega} = (2z - 3xz)dx + (2x^{2} -x - z )dz $$
 with dual vector field
$$ X = (2x^{2} -x - z )\frac{\partial }{\partial x} + (- 2z + 3xz)\frac{\partial }{\partial z}. $$
The singular set of $\G$ is given by
$$ \mathrm{ Sing}( X) =  \Big\{ p_{1} = (2/3, 1, 2/9) ; p_{2} = (1/2, 1, 0) ; p_{3} = (0, 1, 0) \Big\}. $$
We know how to calculate the Grothendieck residue of the foliation $\G$:
$$
 \mathrm{Res}\displaystyle _{c_{1}^{2}}(\G; p_{1}) = \frac{c_{1}^{2}(JX(p_{1}))}{\det(JX(p_{1}))} = \frac{25}{6},
$$
$$
\mathrm{Res}\displaystyle _{c_{1}^{2}}(\G; p_{2}) = \frac{c_{1}^{2}(JX(p_{2}))}{\det(JX(p_{2}))} = - \frac{1}{2},
$$
$$
 \mathrm{Res}\displaystyle _{c_{1}^{2}}(\G; p_{3}) = \frac{c_{1}^{2}(JX(p_{3}))}{\det(JX(p_{3}))} =  \frac{9}{2}.
$$
 Now, we will verify the formula  
$$
\displaystyle  c_{1}^{2}(N_{\sF})  =   \mathrm{Res}(\sF, c_{1}^{2} ; \Gamma)[\Gamma]+ \mathrm{Res}(\sF, c_{1}^{2} ; Q)[Q] + \mathrm{Res}(\sF, c_{1}^{2} ; L)[L]
$$
On the one hand,  Since $\det(N_{\sF})=\mathcal{O}_{\mathbb{P}^3}(4)$, then 
$$c_{1}^{2}(N_{\sF}) = c_{1}^{2}(\det(N_{\sF}))=16h^2,$$
  where $h$ represents the hyperplane class.
On the other hand, by the above calculations  and using that  $[\Gamma]=3h^2$ ,   $[Q]=2h^2$ and  $[L]=h$
we have 
$$\sum_{i = 1}^{3} \mathrm{Res}(\sF, c_{1}^{2} ; Z_{i})[Z_{i}].=  \frac{25}{6}[\Gamma] - \frac{ 1}{2}[Q] + \frac{9}{2}[L]= \frac{25}{6}[3h^2] - \frac{1}{2}[2h^2] + \frac{9}{2}[L]=16h^2. $$

\end{example}

\begin{example}
Let $f:M   \dashrightarrow N$ be a dominant  meromorphic map such that $\dim(N)= k+1$ and $\G$ is an one-dimensional foliation on $N$ with isolated singular set $\mathrm{Sing}(\G)$.
Suppose that $f:M  \dashrightarrow N$ is a submersion outside its  indeterminacy locus $Ind(f)$. Then, the induced foliation $\sF=f^*\G$ on  $M$ has codimension $k$ and $\mathrm{Sing}(\sF)=f^{-1} (\mathrm{Sing}(\G))\cup Ind(f)$.
If  $Ind(f)$ has codimension $\geq k+1$, we conclude that  $\mathrm{cod} (\mathrm{Sing}(\sF))\geq k+1$. If $q\in f^{-1}(p)\subset M$ is a regular point of the map $ f:M   \dashrightarrow N$, then 
$$ \mathrm{Res}(f^*\G, \varphi ; f^{-1}(p)) = \mathrm{Res}_{\varphi}(\G; p)[f^{-1}(p)],  $$
where $ \mathrm{Res}_{\varphi}(\G; p)$ represents  the Grothendieck residue at $p\in \mathrm{Sing}(\G).$
In fact,  there exist  open sets  $U\subset M $ and $V \subset N$,  with $q\in f^{-1}(p)\subset U$ and $p\in V$,  such that $U\simeq f^{-1}(p) \times V $. Now, if we take a $(k+1)$-ball $B$ in $V$ then by 
theorem \ref{principal} we have 
$$ \mathrm{Res}(f^*\G, \varphi ; f^{-1}(p)) = \mathrm{Res}_{\varphi}(\G|_{B}; p)[f^{-1}(p)]=\mathrm{Res}_{\varphi}(\G|; p)[f^{-1}(p)]. $$
For instance,  if $f:\mathbb{P}^n   \dashrightarrow ( \mathbb{P}^{k+1}, \G) $ is a rational linear projection and $\G$  is an one-dimensional foliation with isolated singularities. Since $Ind(f)= \mathbb{P}^{k+1}$, then $ \mathrm{cod} (\mathrm{Sing}(f*\G))=k+1$. Therefore 
$$
 \mathrm{Res}(f^*\G, \varphi ;  f^{-1}(p)) = \mathrm{Res}(f^*\G, \varphi ;  \mathbb{P}^{k+1})=  \mathrm{Res}_{\varphi}(\G; p)[ \mathbb{P}^{k+1}].
$$
\end{example}

\section{Cenkl  algorithm for  singularities with non-expected dimension  }

In \cite{Ce}  Cenkl  provided an algorithm to determinate  residues for non-expected dimensional singularities,  under a certain regularity condition on the singular set of the foliation.
We observe that this condition is not necessary. In fact, Cenkl 's conditions are the following:
 
 Suppose that the  singular set $S:=\mathrm{Sing}(\sF)$ of $\sF$ has pure codimension $k+s$, with $s\geq 1$,  and
 \begin{enumerate}
\item[(i)] $\mathrm{cod}(S) \geq 4$.
\item[(ii)] there exists a closed subset  $W \subset M$  such that $S\subset W$ with  the property
$$
\H^j(W, \mathbb{Z}) \simeq \H^j(W\setminus S , \mathbb{Z}), \  \  j=1,2.  
$$
 \end{enumerate}
Denote by $M'=M\setminus S$, Cenkl show that under the above  condition the line bundle $\wedge^k( N_{\sF}|_{M'}^{\vee} )$  on $M'$ can be extended a line bundle on $M$.
We observe that there always exists a line bundle $\det(N_{\sF})^{\vee }= [\wedge^k (N_{\sF})^\vee ]^{\vee \vee}$ on $M$ which extends  $\wedge^k( N_{\sF}|_{M'}^{\vee} )$, since $N_{\sF}$ is a torsion free sheaf and $S=\mathrm{Sing}(N_{\sF})$. 
See,  for example  \cite[Proposition 5.6.10 and Proposition 5.6.12  ]{Ko}.
Now, consider the vector bundle
$$
E_{\sF}= \det(N_{\sF})^{\vee }\oplus \det(N_{\sF})^{\vee }.
$$ 
Observe that $E_{\sF}|_{M'}= \wedge^k( N_{\sF}|_{M'}^{\vee} ) \oplus  \wedge^k( N_{\sF}|_{M'}^{\vee} ) $.  
Thus, we conclude that  Lemma 1 in  \cite{Ce} holds in general:
\begin{lemma}\label{lemaCe}
Consider the projective bundle $\pi: \mathbb{P}(E_{\sF}) \to M$.
 Then  there exist a holomorphic foliation $\sF_{\pi}$ on $ \mathbb{P}(E_{\sF}) $  with   singular set  $\mathrm{Sing}(\sF_{\pi})=\pi^{-1}(S)$ such that 
\begin{center}
$\dim(\sF_{\pi})=\dim(\sF)$ and $\dim(\mathrm{Sing}(\sF_\pi))=\dim(S)+1.$
\end{center}
\end{lemma}

We succeeded in replacing the compact manifold M with a foliation  $\sF$ and the singular set $S$ such that 
$\dim(\sF)- \dim(\mathrm{Sing}(\sF))=n-s$ by another compact manifold $\mathbb{P}(E_{\sF})$ and a foliation $\sF_{\pi}$ with singular set 
$\dim(\sF_\pi)- \dim(\mathrm{Sing}(\sF_\pi))=n-s-1$
. If this procedure is repeated $(n-s-1)$-times we end up with a compact complex analytic manifold with a holomorphic foliation whose singular set is a subvariety of complex dimension one less than the leaf dimension of the foliation. That is,  we have a tower of foliated manifolds 
\begin{equation}\nonumber
\xymatrix{
 (P_{n-s-1}, \sF^{n-s-1})    \ar[r]^{\pi_{n-s-1}}    &  (P_{n-s-2}, \sF^{n-s-2})  \ar[r]            & \cdots \ar[r]^{\pi_{2}}     & (P_{1}, \sF^{1})   \ar[r]^{\pi_1:=\pi}     &    (M,\sF)
  }
\end{equation}
where $(P_{i}, \sF^{i}) $ is such that $P_{i}=  \mathbb{P}(E_{\sF^{i-1}})$ and $(P_{1}, \sF^{1})=( \mathbb{P}(E_{\sF}), \sF_{\pi}) $. Thus, by Lemma \ref{lemaCe} we conclude that on $ P_{n-s-1}$ we have a foliations $\sF^{n-s-1}$ such that $\mathrm{Sing}(\sF^{n-s-1})=(\pi_{n-s-1} \circ \cdots \circ \pi_2\circ\pi_1)^{-1}(S) $ and 
$$
\dim(\mathrm{Sing}(\sF^{n-s-1}))=\dim(\sF^{n-s-1})-1.
$$
That is, $\mathrm{cod}(\mathrm{Sing}(\sF^{n-s-1}))= \mathrm{cod}(\sF^{n-s-1})+1$. 

On the one hand, we can apply the Theorem \ref{principal} to determinate the residues of $\sF^{n-s-1}$. On the other hand, Cenkl show that we can 
calculate the residue $\mathrm{Res}_{\varphi}(\sF^1, Z_1)$ in terms of the residue $\mathrm{Res}_{\varphi}(\sF, Z)$  for symmetric polynomial $\varphi$ of degree $k+1$. 

Let us recall the Cenkl's construction:

 Let $\sigma_1,\dots,\sigma_{\ell}$ be the elementary symmetric functions in the $n$ variables $x_1,\dots,x_n$ and let $\rho_1,\dots, \rho_{\ell}$   be the elementary symmetric functions 
 in the $n+1$ variables  $x_1,\dots,x_n,y$. It follows from \cite[Corollary, pg 21]{Ce} that for   any polynomial  $\phi$, of degree $\ell$, can be associated a polynomial $\psi$ of degree $\ell +1$ such that
 $$
\psi(\rho_1,\dots, \rho_{\ell} )  =   \phi(\sigma_1,\dots,\sigma_{\ell}) y + \phi^0(\sigma_1,\dots,\sigma_{\ell}) + \sum_{j\geq 2} \phi^j(\sigma_1,\dots,\sigma_{\ell}) \cdot y^j,
 $$
 where $\phi^0$ has degree $\ell +1$  and $ \phi^j$ has degree $\ell -j+1$.

Let $T_{P/M}$ be  the tangent bundle associated the one-dimensional foliation induced by the $\mathbb{P}^1$-fibration   $(P, \sF_{\pi})  \to     (M,\sF)$.

 Therefore,  it follows from Lemma \ref{lemaCe}, Cenkl's construction \cite[Theorem 1]{Ce} and Theorem \ref{principal}    the following : 

\begin{thm}
Suppose that $\mathrm{cod}(\mathrm{Sing}(\sF))\geq \mathrm{cod}(\sF )+2$.   If  $\varphi$ is  a homogeneous symmetric polynomials
of degree $\mathrm{cod}(\sF )+1$,  then
$$
\mathrm{Res}_{\psi}(\sF^1|_{B_p}; p)[Z_1]   =  \pi^*\mathrm{Res}_{\varphi}(\sF, Z) \cap c_1(T_{P/M}) + \pi^*(\phi^0(N_{\sF})) + \sum_{j\geq 2}  \pi^*(\phi^j(N_{\sF}))\cap c_1(T_{P/M})^j,
$$
 where  $\mathrm{Res}_{\psi}(\sF^1|_{B_p}; p)$ represents  the Grothendieck residue at $p$ of the one dimensional  foliation $\sF^1|_{B_p}$ on a $(k+1)$-dimensional  transversal ball $B_p$.
\end{thm}

We believe  that  this algorithm can be adapted to  the context of   residues for flags of foliations \cite{BCL}.


\begin{thebibliography}{999}

\bibitem{BB1} P. Baum, R. Bott, {\it On the zeros of meromorphic vector fields, Essay on Topology and Related Topics}, Spring-Verlag, New York, 1970, 29-47.


\bibitem {BB}   P. Baum, R. Bott,  \newblock
{\it Singularities of holomorphic foliations}, J. Differential Geom. 7 (1972), 279-342.

\bibitem {BS} 
F. Bracci, T. Suwa, { \it Perturbation of Baum-Bott residues}. Asian J. Math., 19, 5, (2015), 871-886

\bibitem {BCL} 
J-P. Brasselet, M. Corr\^ea, F. Louren\c co, { \it Residues for  flags of holomorphic foliations}.  Advances in Mathematics, vol. 320, n.7, 1158-1184,  2017.

\bibitem {BruPer}  M.  Brunella and C. Perrone, 
{\it Exceptional singularities of codimension one holomorphic
foliations}, Publicacions Matem\`atiques 55 (2011), 295-312.

\bibitem {CerLins}  D.  Cerveau, A.  Lins Neto, 
{\it Irreducible components of the space of holomorphic foliations of degree two in  CP(n) }, Ann. Math. 143 (1996) 577-612.

\bibitem {Ce} 
B. Cenkl, { \it Residues of singularities of holomorphic foliations}, J. of Differential
Geometry, 13 (1978) 11-23.

 



\bibitem {CoFer}  M. Corr\^ea,  A. Fernand\'ez-P\'erez,
{\it Absolutely $k$-convex domains and holomorphic foliations on homogeneous manifolds},  Journal of the Mathematical Society of Japan. vol. 69, n.3,  1235-1246, 2017.


\bibitem{gri} P. Griffiths, J. Harris,  {\it Principles of algebraic geometry}, Wiley, 1978.


\bibitem {Jou} J-P.  Jouanolou,  \newblock
{\it Equations de Pfaff alg\'ebriques}, 1979 Lecture Notes in Mathematics, 708. Springer-Verlag, Berlin.

\bibitem {Ko}
S. Kobayashi, {\it Differential geometry of complex vector bundles}, Princeton
Univ. Press, 1987.


\bibitem {Soa}M. Soares, 
{\it Holomorphic foliations and characteristic numbers}, Comm. Contemporary Maths. 7(5) (2005), 583-596.



\bibitem {Suwa}   T. Suwa,  
{\it Indices of Vector Fields and Residues of Singular Holomorphic
Foliations}, Actualit\'es Math\'ematiques, Hermann, Paris 1998.

\bibitem {Suwa1} T.  Suwa,  
{\it Residues of Complex analytic Foliations Singularities}, J. Math. Soc. Japan., 36 (1984), 37-45.



\bibitem {Vis}  M. S. Vishik, 
 {\it Singularities of analytic foliations and characteristic classes}, Functional Anal. Appl. 7 (1973) 1-15.


\end{thebibliography}
\end{document}